\documentclass[12pt]{amsart}
\usepackage{amsmath,amssymb,amsbsy,amsfonts,amsthm,latexsym,
            amsopn,amstext,amsxtra,euscript,amscd,amsthm,graphicx,comment}
\usepackage[left=3.5cm,right=3.5cm,top=2cm,bottom=2cm,includeheadfoot]{geometry}
\usepackage[enableskew]{youngtab}

\newtheorem{thm}{Theorem}
\newtheorem{theorem}[thm]{Theorem}

\newtheorem{cor}{Corollary}
\newtheorem{corollary}[cor]{Corollary}

\newtheorem{conj}{Conjecture}
\newtheorem{conjecture}[conj]{Conjecture}

\newtheorem{defi}{Definition}
\newtheorem{definition}[defi]{Definition}

\theoremstyle{remark}
\newtheorem{remark}{Remark}

\def\\{\cr}
\def\({\left(}
\def\){\right)}
\def\<{\langle}
\def\>{\rangle}

\begin{document}
\title[Conjectures Nekrasov--Okounkov hook length formula]
{On conjectures regarding the Nekrasov--Okounkov hook length formula}
\author{Bernhard Heim }
\address{German University of Technology in Oman (GUtech), 
PO Box 1816, Athaibah PC 130
, Sultanate of Oman}
\address{Faculty of Mathematics, Computer Science, and Natural Sciences,
RWTH Aachen University, 52056 Aachen, Germany}
\email{bernhard.heim@gutech.edu.om}
\author{Markus Neuhauser}
\address{German University of Technology in Oman (GUtech),
PO Box 1816, Athaibah PC 130
, Sultanate of Oman}
\address{Faculty of Mathematics, Computer Science, and Natural Sciences,
RWTH Aachen University, 52056 Aachen, Germany}
\email{markus.neuhauser@gutech.edu.om}
\pagenumbering{arabic}

\begin{abstract}
The Nekrasov--Okounkov hook length formula provides a fundamental link between the theory of partitions and the
coefficients of powers of the Dedekind eta function. 
In this paper we examine three 
conjectures presented by Amdeberhan.
The first conjecture is a refined Nekrasov--Okounkov formula involving hooks with trivial legs. We prove the conjecture.
The second conjecture is on properties of the roots of the underlying D'Arcais polynomials. 
We give a counterexample and present a new conjecture.
The third conjecture is on the unimodality of the coefficients of the involved polynomials. 
We confirm the conjecture up to the polynomial degree $1000$.
\end{abstract}
\maketitle
\section{Introduction}
In 2006 Nekrasov and Okounkov \cite{NO06, We06,Ha10} published a new type of hook length formula.
Based on their work on random partitions and the Seiberg--Witten theory they 
obtained an unexpected identity relating the sum over products of partition hook lengths \cite{Ma95,Fu97} to 
the coefficients of complex powers of Euler products \cite{Ne55,Se85,HNW18}, which is essentially a power of the Dedekind eta function.
This paper is devoted to three open conjectures stated by Amdeberhan in \cite{Am15} (section 2).

Let $\lambda$ be a partition of $n$ denoted by $\lambda \vdash n$ with weight $|\lambda|=n$. We denote by
$\mathcal{H}(\lambda)$ the multiset of hook lengths associated to $\lambda$ and $\mathcal{P}$ be the set of all partitions.
The Nekrasov--Okounkov hook length formula is given by
\begin{equation} \label{ON}
\sum_{ \lambda \in \mathcal{P}} q^{|\lambda|} \prod_{ h \in \mathcal{H}(\lambda)} \left(  1 - \frac{z}{h^2} \right) =   \prod_{m=1}^{\infty} \left( 1 - q^m \right)^{z-1}.
\end{equation}
Let $q:=e^{2\pi i\tau }$ and $ \tau$ be in the upper complex half-plane. The identity (\ref{ON}) is valid for all $z \in \mathbb{C}$.
The Dedekind eta function $\eta(\tau)$ is given by $q^{\frac{1}{24}} \prod_{m=1}^{\infty} \left( 1 - q^m \right)$ (see \cite{On03}).
\newline
\newline
The first conjecture is a refinement of the Nekrasov--Okounkov hook length formula \cite{Am15}.
\begin{conjecture} \label{conjecture1}
Let $\mathcal{H}(\lambda)^{\diamond}$ be the
multiset of hook lengths with trivial legs. Then
\begin{equation} \label{mainexpression}
\sum_{ \lambda \vdash n} 
\prod_{ h \in \mathcal{H}(\lambda)^{\diamond}} 
\left(  \frac{h + z}{h} \right) 
=
\sum_{ \lambda \vdash n} 
\prod_{ h \in \mathcal{H}(\lambda)} 
\left(  \frac{h^2 + z}{h^2} \right). 
\end{equation}
\end{conjecture}

The second conjecture is on the roots of the polynomials given in (\ref{mainexpression}).
Note that this is also related to the Lehmer conjecture \cite{Le47, HNW18}.
\begin{conjecture} \label{conjecture2}
Let $n$ be a positive integer. Then the polynomial
\begin{equation}\label{Q}
Q_n(z) := 
\sum_{ \lambda \vdash n} 
\prod_{ h \in \mathcal{H}(\lambda)} 
\left(  \frac{h^2 + z}{h^2} \right) \in \mathbb{C}[z]
\end{equation}
has (i) only simple roots, (ii) only real roots, (iii) only negative roots.
\end{conjecture}

The third conjecture is on the coefficients of the polynomials defined in (\ref{Q}).
\begin{conjecture} \label{conjecture3}
Let $n$ be a positive integer. Then $Q_n(z)$ is unimodal.
\end{conjecture}
In 1913 D'Arcais \cite{DA13} studied a sequence of polynomials $P_n(x)$:
\begin{equation}
\sum_{n=0}^{\infty} P_n(z) \, q^{n}
=  \prod_{n=1
}^{\infty} \left( 1 - q^n \right)^{-z}. \label{Arcais}
\end{equation}
The coefficients are called D'Arcais numbers \cite{Co74}.
Independently from D'Arcais, Newman and Serre \cite{Ne55, Se85} studied the polynomials in the context of modular forms.
Serre proved a famous theorem on lacunary modular forms, utilizing the factorization of $P_n(x)$ for $ 1 \leq n \leq 10$ over $\mathbb{Q}$.
All authors mentioned so far introduced slightly different normalized polynomials: D'Arcais \cite{DA13},  Newman and Serre \cite{Ne55, Se85}, Amdeberhan \cite{Am15}.
D'Arcais definition fits best to the new Conjecture 2. 

In this paper we prove Conjecture 1. We use the Nekrasov--Okounkov hook length formula and
combinatorial arguments. It is sufficient to show the formula evaluated at negative integer values.
We give a 
counterexample to Conjecture 2. In the degree $10$ case, $P_{10}(x)$, non-real complex roots exist.
The simplicity of the roots was already studied in \cite{HN18}. Finally we give numerically evidence for Conjecture 3.

\section{New partition Hook length formula}
A partition of $n$ (for an introduction we refer to \cite{Ma95,St99,AE04,Ha09})
is a finite decreasing sequence 
$\lambda = (\lambda_1,\lambda_2, \ldots, \lambda_l)$ of positive integers
such that $\vert \lambda \vert := \sum_j \lambda_j =n$. We write $\lambda \vdash n$.
The set of all partitions of $n$ is denoted $\mathcal{P}(n)$ and the set of all partitions for all $n \in \mathbb{N}$ is
denoted $\mathcal{P}$.

The integers $\lambda_j$ are called the parts of $\lambda$ and $l= l(\lambda)$ the length of the partition.
Partitions are presented by their Young diagram.
Let $\lambda = (7,3,2)$. Then $l
\left( \lambda\right) = 3$ and $ n = \vert \lambda \vert = 12$.
$$
\young(~~~~~~~,4~~,~~). 
$$
We attach to each cell $u$ of the diagram the arm $a_u(\lambda)$, the amount of cells in the same row of $u$ to the right of $u$.
Further we have the leg $\ell _{u}\left( \lambda\right) $, the number of cells in the same column of $u$ below of $u$.
The hook length $h_u(\lambda)$ of the cell $u$ is given by $h_u(\lambda):= a_u(\lambda) + \ell _{u}\left( \lambda\right) + 1$.
The hook length multiset $\mathcal{H}(\lambda)$ is the multiset of all hook length of $\lambda$.
Our example gives 
$$\mathcal{H}(\lambda) = \{ 9,8,6,4,3,2,1,4,3,1,2,1 \}.$$

The list is 
given from the left to the right and from the top to the bottom in the Young diagram.
Cells have the coordinates $(i,j)$ following the same procedure, hence $4$ is in the cell $(2,1)$ with $\ell =1$ and $a=2$
(we usually simplify the notation).
Let $f_{\lambda}$ denote the number of standard Young tableaux of shape $\lambda$. 
These are all possible combinations of filling a Young diagram 
with the numbers $\{1,2,\ldots,n\}$ for $\lambda \vdash n$, such
that each number occurs once and in each row and column 
(from the left to the right and from the top to the bottom) the numbers are strictly increasing. The famous classical hook formula (Frame, Robinson, Thrall) states:
\begin{equation}
f_{\lambda} = \frac{n!}{\prod_{h \in \mathcal{H}(\lambda)} \, h}.
\end{equation}
In the following we prove Conjecture 1.
Let $\mathcal{H}(\lambda)^{\diamond}$ be the multiset of hook lengths with trivial legs. Then
Conjecture 1 states:
\begin{equation} 
Q_n(z):= \sum_{ \lambda \vdash n} 
\prod_{ h \in \mathcal{H}(\lambda)^{\diamond}} 
\left(  \frac{h + z}{h} \right) 
= P_n(z+1).
\end{equation}
It can easily be verified for $z=-1$ and $z=0$.
Since there always  exists 
a hook length $h=1$ for $\lambda$. The product in $Q_n(-1)$ always has a factor
which is zero. Hence $Q_n(-1)= P_n(0)=0$ for $n \geq 1$.
Let $z=0$. Then $$Q_n(0) = \sum_{ \lambda \vdash n} 1,$$ 
which is equal to the number of partitions of $n$. By Euler it is known that this is the coefficient of $q^n$
in the product $\prod_{m=1}^{\infty }\left( 1-q^m\right) ^{-1}$, hence equal to $P_n(1)$.
\begin{proof}
[Proof of Conjecture 1.]
\
\newline
The Nekrasov--Okounkov hook length formula states that 
\begin{equation} 
\sum_{ \lambda \vdash n} 
\prod_{ h \in \mathcal{H}(\lambda)} 
\left(  \frac{h^2 + z}{h^2} \right) 
\end{equation}
is equal to $P_n(z+1)$. Hence it is sufficient to prove that
\begin{equation}
\sum_{ \lambda \vdash n} 
\prod_{ h \in \mathcal{H}(\lambda)^{\diamond}} 
\left(  \frac{h -(m+1)}{h} \right) 
= P_n(-m)
\end{equation}
for all $m \in \mathbb{N}$.
Let $\lambda = ( \lambda_1, \ldots, \lambda_l)$ be a partition of $n$.
We count all parts of $\lambda$ with value $j$ and put:
$$a_j (\lambda):= \sharp \{ i \, | \, \lambda_i = j \}.$$ This leads to the bijection
\begin{equation}
 \psi: \big\{ \lambda \vdash n \big\} \longrightarrow \left\{ (a_1, \ldots, a_n) \in \mathbb{N}_0 \, | \,\sum_{k=1}^n k \, a_k = n \right\},
\end{equation}
where $\lambda$ maps to $a(\lambda) = \left( a_1(\lambda), \ldots, a_n(\lambda) \right) $.
We collect all the terms in $$\prod_{d=1}^{\infty} ( 1 - q^d)^m$$ contributing to $P_n(-m)$. Then $\psi(\lambda)$
contributes with multiplicity 
$$  \binom{m}{a_1(\lambda)} \cdots 
\binom{m}{a_n(\lambda)}.$$
Hence we obtain
\begin{equation} \label{first} P_n(-m) = \sum_{ \lambda \vdash n} (-1)^{
l\left( \lambda \right) } \,\,  \binom{m}{a_1(\lambda)} \cdots 
\binom{m}{a_n(\lambda)}. \end{equation}
Next we study the hook length term
\begin{equation}
C_m(n) := \sum_{ \lambda \vdash n} 
\prod_{ h \in \mathcal{H}(\lambda)^{\diamond}} 
\left(  \frac{h -(m+1)}{h} \right). 
\end{equation}
Note that if  $h \in \mathcal{H}(\lambda)^{\diamond}$, then also $1,\ldots, h-1$ are elements of $\mathcal{H}(\lambda)^{\diamond}$ for $h > 1$.
Let $\lambda^c$ denote the conjugate of $\lambda$, 
which is also a partition of $n$. Since $\{ \lambda \vdash n \} = \{ \lambda^c \, \, \vert \, \lambda \vdash n \}$, we have
\begin{equation}
 \sum_{ \lambda \vdash n} 
\prod_{ h \in \mathcal{H}(\lambda)^{\diamond}} 
\left(  \frac{h -(m+1)}{h} \right) =   \sum_{ \lambda^c \vdash n} 
\prod_{ h \in \mathcal{H}(\lambda)^{\diamond}} 
\left(  \frac{h -(m+1)}{h} \right).
\end{equation}
We consider the Young diagram of $\lambda^c$. For example let $\lambda = (4,3,3,2
,1,1)$. Then we have
$$\young(~~~~21,~~~1,~21,1).
$$
The numbers denote the hook length for all cells with no leg. 
Let $$b_j(\lambda) :=  \sharp \{ i \, | \, \lambda_i = j\}.$$
Then $b_1(\lambda), \ldots, b_n(\lambda) $ contributes to $C_m(n)$ with multiplicity
$$  \binom{m}{b_1(\lambda)} \cdots 
\binom{m}{b_n(\lambda)}.$$
Here we used the simple identity $$\prod_{k=1}^h  \left( \frac{ k
- (m+1)}{k
} \right) = (-1)^h \binom{m}{h}.$$
Hence we obtain 
\begin{equation}\label{second}
C_m(n) = \sum_{ \lambda \vdash n} (-1)^{l(\lambda)} \,\,
\binom{m}{b_1(\lambda)} \cdots 
\binom{m}{b_n(\lambda)}.
\end{equation}
Comparing (\ref{first}) and (\ref{second}) proves Conjecture 1.
\end{proof}
\begin{corollary}
Let $\mathcal{H}(\lambda)^{\diamond \diamond}$ be the multiset of hook lengths with trivial arms. Then
\begin{equation} 
\sum_{ \lambda \vdash n} 
\prod_{ h \in \mathcal{H}(\lambda)^{\diamond}} 
\left(  \frac{h + z}{h} \right) 
=
\sum_{ \lambda \vdash n} 
\prod_{ h \in \mathcal{H}(\lambda)^{\diamond \diamond}}
\left(  \frac{h + z}{h} \right).
\end{equation}
\end{corollary}
Let $ \lambda \vdash n$. Let $k_j:= k_j(\lambda) = a_j(\lambda)= b_j(\lambda)$.
Since 
\begin{equation}
\binom{k+z}{k} = (-1)^k \binom{-z-1}{k}
\end{equation}
for $k \in \mathbb{N}_0$ and $z \in \mathbb{C}$, we obtain (put $z=-(m+1)$):
\begin{corollary}
Let $\mathcal{H}(\lambda)^{\diamond}$ be the multiset of hook lengths with trivial legs. Then
\begin{equation*} 
\sum_{ \lambda \vdash n} 
\prod_{ h \in \mathcal{H}(\lambda)^{\diamond}} 
\left(  \frac{h + z}{h} \right) 
=
\sum_{ \lambda \vdash n} 
\prod_{j=1}^{n}
\binom{k_j +z }{k_j} =
\sum_{ \lambda \vdash n} 
\prod_{j=1}^{n} (-1)^{
l\left( \lambda \right) }
\binom{-z-1}{k_j}.
\end{equation*}
\end{corollary}
This proves the complete Conjecture 2.1 given in \cite{Am15}.
It would be also interesting to study the Conjecture  in terms of $t$-cores (\cite{GKS90}).
\section{D'Arcais type polynomials}
The Nekrasov--Okounkov hook length formula \cite{NO06}, \cite{Ha10} states:
\begin{equation}
\sum_{ \lambda \in \mathcal{P}} q^{|\lambda|} \prod_{ h \in \mathcal{H}(\lambda)} \left(  1 - \frac{z}{h^2} \right) =  \sum_{n=0}^{\infty} P_n(1-z) \, q^n.
\end{equation}
The polynomials $P_n(x)$ can also be recursively defined by $$P_n(x):= \frac{x}{n} \sum_{k=1}^n \sigma(k) P_{n-k}(x).$$ Here $P_0(x):=1$ and $\sigma(n):= \sum_{d \mid 
n} d$.
This makes it possible to calculate the coefficients of the polynomials directly. The first $20$ polynomials can be found in \cite{HNR18}. 
Note that the first $10$ were published in 1955 by Newman \cite{Ne55} (see also \cite{Se85}).
We claim that $P_{10}(x)$ has non-real roots. This was tested numerically in \cite{HNR18}.
In this paper we use a Theorem of Aissen--Schoenberg--Whitney and Edrai \cite{ASW52} to give a rigorous algebraic proof. 
This leads to a counterexample for Conjecture 2.
We also give a second proof involving derivatives.

Let us first recall the definition of a totally nonnegative matrix. Then we give a bijection between the set of
infinite real sequences and certain Toeplitz matrices.  We define (finite) Polya frequencies sequences and state a
well known characterization of polynomials with real roots.

Let $A=(a_{i,j=0}^{\infty})$ be a two way infinite matrix with real entries. 
The matrix $A$ is called totally nonnegative (TN) if all minors are nonnegative.

\begin{definition}
Let $a_0,a_1,\ldots $ be an infinite sequence of real numbers. We attach the matrix 
$A=(a_{i,j})_{i,j=0}^{\infty}$ given by $a_{i,j}:= a_{i-j}$ for $0 \leq j \leq i $ and $a_{i,j}=0$ otherwise.
\end{definition}
\begin{
example}
Let the sequence $2,2,1,0,0, \ldots $ be given, then the attached matrix $A$ is given by
$$
\begin{pmatrix}
2 & 0& 0 &0 & 0& \dots \\
2 & 2& 0 &0 & 0& \dots  \\
1 & 2& 2 &0 & 0& \dots   \\
0 & 1& 2 &2 & 0&  
\\
0 & 0& 1 &2 & 2& 
\\
\vdots & \vdots & 
&  \ddots & \ddots & \ddots   
\end{pmatrix}.
$$
\end{example}
We consider finite sequences as infinite sequences. Note that $A$ is not a TN, since the minor, the determinant of
$$
\begin{pmatrix}
2 & 2& 0 &0 \\
1 & 2& 2 &0 \\
0 & 1& 2 &2 \\
0 & 0& 1 &2  
\end{pmatrix},
$$
is negative.
\begin{definition}
Let $a_0,a_1,\ldots ,a_n$ be a finite sequence of nonnegative real numbers. 
Then the sequence is called Polya frequency sequence if the attached (infinite) matrix $A$ is a totally nonnegative matrix.
\end{definition}
The following result is well known \cite{KRY09}. It gives a sufficient and necessary criterion for polynomials with nonnegative coefficients to have only real roots.
Note that this can be applied to the D'Arcais polynomials $P_n(x)$.
We recall \cite{HLN18} that $P_n(x) = \frac{x}{n!} \sum_{k=0}^{n-1} a_k^{(n)} \, x^k$, where $a_k^{(n)} \in \mathbb{N}$. Here $a_{n-1}^{(n)}=1$.
\begin{theorem}[Aissen--Schoenberg--Whitney, Edrai]
A finite sequence $a_0,\ldots , a_n$ of nonnegative real numbers is a Polya frequency sequence if and only if the attached polynomial $\sum_{i=0}^n a_i x^i$ has only real roots.
\end{theorem}
We have shown that the example $(2,2,1)$ is not a Polya frequency sequence, which directly reflects that the polynomial $2 + 2x + x^2$ has a non-real root.

Now we apply the theorem to the polynomial $P_{10}(x)$. Actually \cite{Ne55,HNR18} $$P_{10}(x) = \frac{x}{10!} \, (x+1) \, R(x).$$ Here $R(x) = \sum_{k=0}^8 a_k \, x^k $ with $a_k \in \mathbb{N}$.
Hence it is sufficient to show that the coefficients of $R(x)$ are not a Polya frequency sequence.
%
\begin{eqnarray*}
a_0= 6531840, \,
a_1 = 29758896,\,
a_2 = 28014804,\,
a_3 = 10035116,\\
a_4 = 1709659, \,
a_5 = 147854,\,
a_6 = 6496, \,
a_7 = 134, \, 
a_8 = 1. \, 
\end{eqnarray*}
Then the attached matrix $A$ (infinite rows and columns) is given by
{\tiny
$$
 \left[ \begin {array}{cccccccccc} 6531840&0&0&0&0&0&0&0&0& \dots
\\ \noalign{\medskip}29758896&6531840&0&0&0&0&0&0&0& \dots
\\ \noalign{\medskip}28014804&29758896&6531840&0&0&0&0&0&0& \dots
\\ \noalign{\medskip}10035116&28014804&29758896&6531840&0&0&0&0&0& \dots
\\ \noalign{\medskip}1709659&10035116&28014804&29758896&6531840&0&0&0&0& \dots
\\ \noalign{\medskip}147854&1709659&10035116&28014804&29758896&6531840
&0&0&0& \dots
\\ \noalign{\medskip}6496&147854&1709659&10035116&28014804&
29758896&6531840&0&0 &\dots\\ 
\noalign{\medskip}134&6496&147854&1709659&
10035116&28014804&29758896&6531840&0&\\ 
\noalign{\medskip}1&134&6496& 
147854&1709659&10035116&28014804&29758896&6531840 & \\
\noalign{\medskip} 0 & 1&134&6496& 
147854&1709659&10035116&28014804&29758896&
\ddots\\
&&\ddots&\ddots&\ddots&\ddots&\ddots&\ddots&\ddots&\ddots
\end {array} \right]
$$}
%
A matrix  that violates the condition of the theorem
(i.~e.\ has negative determinant) is the
following $26\times 26$ matrix:
{\tiny
\[
\left[
\begin{array}{cccccccccc}
\\ \noalign{\medskip}10035116&28014804&29758896&6531840&0&0&0&0&0& \dots
\\ \noalign{\medskip}1709659&10035116&28014804&29758896&6531840&0&0&0&0& \dots
\\ \noalign{\medskip}147854&1709659&10035116&28014804&29758896&6531840
&0&0&0& \dots
\\ \noalign{\medskip}6496&147854&1709659&10035116&28014804&
29758896&6531840&0&0 &\dots\\ 
\noalign{\medskip}134&6496&147854&1709659&
10035116&28014804&29758896&6531840&0&
\\ 
\noalign{\medskip}1&134&6496& 
147854&1709659&10035116&28014804&29758896&6531840 & \\
\noalign{\medskip} 0 & 1&134&6496& 
147854&1709659&10035116&28014804&29758896&
\ddots\\
\noalign{\medskip} 0 & 0 & 1&134&6496& 
147854&1709659&10035116&28014804&
\ddots\\
\noalign{\medskip} 0 & 0 & 0 & 1&134&6496& 
147854&1709659&10035116&
\ddots\\
\vdots&\vdots&&&\ddots&\ddots&\ddots&\ddots&\ddots&\ddots
\end{array}
\right]
\]
}
It contains the rows~$4$ to $29$
and columns~$1$ to $26$ of the
infinite 
Toeplitz matrix $A$. Hence $P_{10}(x)$ has non-real roots.

For the convenience of the reader we give a second proof using the derivative $R'(x)$.
Inserting the limiting values of the following intervals
into $R\left( x\right) $ we observe
that $R\left( x\right) $ has (at least)
one root 
$z_{1}\in \left( -59,-58\right) $,
$z_{2}\in \left( -33,-32\right) $,
$z_{3}\in \left( -18,-17\right) $,
$z_{4}\in \left( -14,-13\right) $,
$z_{7}\in \left( -2,-1\right) $, and
$z_{8}\in \left( -1,0\right) $.
Similarly for the derivative
$R^{\prime }\left( x\right) $ we find
a root
$z_{1}^{\prime }\in \left( -53,-52\right) $,
$z_{2}^{\prime }\in \left( -29,-28\right) $,
$z_{3}^{\prime }\in \left( -16,-15\right) $,
$z_{4}^{\prime }\in \left( -11,-10\right) $,
$z_{5}^{\prime }\in \left( -6,-5\right) $,
$z_{6}^{\prime }\in \left( -4,-3\right) $, and
$z_{7}^{\prime }\in \left( -1,0\right) $.
The degree of $R^{\prime }\left( x\right) $ is~$7$.
Hence
each of these intervals contains exactly one
(simple) root.

Firstly we note that
in particular the root
$z_{7}^{\prime }\in \left( -1,0\right) $ of
$R^{\prime }\left( x\right) $ is unique,
$R^{\prime }\left( x\right) $ does not have a root
in $\left( -2,-1\right] $,
and there is a root $z_{7}\in \left( -2,-1\right) $
of $R\left( x\right) $.
Since the roots of
$R^{\prime }\left( x\right) $
are simple
$z_{7}^{\prime }$ could be
a double root of
$R\left( x\right) $. But
this would contradict
the opposite signs of the
values of
$R\left( x\right) $ for the
limits of
$\left( -1,0\right) $.
(Note that there is only
\emph{one\/} root of
$R^{\prime }\left( x\right) $
in $\left( -1,0\right) $.)
Hence we must have $z_{7}<z_{7}^{\prime }<z_{8}$.
Thus a root smaller than $z_{1}$ and larger than $z_{8}$ of
$R\left( x\right) $ is not possible since this
would imply a root of $R^{\prime }\left( x\right) $
for $x<z_{1}$ or $x>z_{8}$ and we have found
that there are no such roots.

From the distribution of the roots of
$R^{\prime }\left( x\right) $ we can observe
that $R\left( x\right) $ is increasing on
$\left( z_{3}^{\prime },z_{4}^{\prime }\right) $, decreasing on
$\left( z_{4}^{\prime },z_{5}^{\prime }\right) $, and increasing on
$\left( z_{5}^{\prime },z_{6}^{\prime }\right) $.
Since $R\left( -6\right) >0$ and $R\left( -5\right) >0$
with a minimum at $z_{5}^{\prime }$ in between the only
chance for the `missing' zeros of $R\left( x\right) $
would be 
to be in the interval $\left( -6,-5\right) $.

We show next that $R\left( x\right) $
is strictly positive on $\left( -6,-5\right) $.
Thus the second derivative
$R^{\prime \prime }\left( x\right) $ has exactly one root
in
$\left( z_{4}^{\prime },z_{5}^{\prime }\right) \supset \left( -10,-6\right) $
and
$\left( z_{5}^{\prime },z_{6}^{\prime }\right) \supset \left( -5,-4\right) $.
The checking of the limiting values shows that there is exactly one
root of $R^{\prime \prime }\left( x\right) $
in $\left( -9,-8\right) $ and $\left( -5,-4\right) $.
This means for the third derivative
$R^{\prime \prime \prime }\left( x\right) $
that in between these two roots of
$R^{\prime \prime }\left( x\right) $
there is exactly one root
$R^{\prime \prime \prime }\left( x\right) $.
Checking the limiting values shows that this root is
in the interval $\left( -7,-6\right) $, which means
that there is no sign change of
$R^{\prime \prime \prime }\left( x\right) $ on
$\left( -6,-5\right) $ i.~e.\ 
$R^{\prime \prime \prime }\left( x\right) <0$
for $x\in \left( -6,-5\right) $.
Expanding $R\left( x\right) $ around $-5$ yields
$R\left( x\right) = 1632960 + 1690056\left( x+5\right) + 1663164\left( x+5\right) ^2+$ terms of higher degree.
Since the third derivative is negative on
$\left( -6,-5\right) $, we obtain
$R\left( x\right) > 1632960 + 1690056\left( x+5\right) + 1663164\left( x+5\right) ^2$
for $x\in \left( -6,-5\right) $.
The discriminant of this quadratic polynomial is
negative so it does not have real roots and the same
holds for $R\left( x\right) $ on
$x\in \left( -6,-5\right) $.

\begin{figure}
\includegraphics[width=12cm]{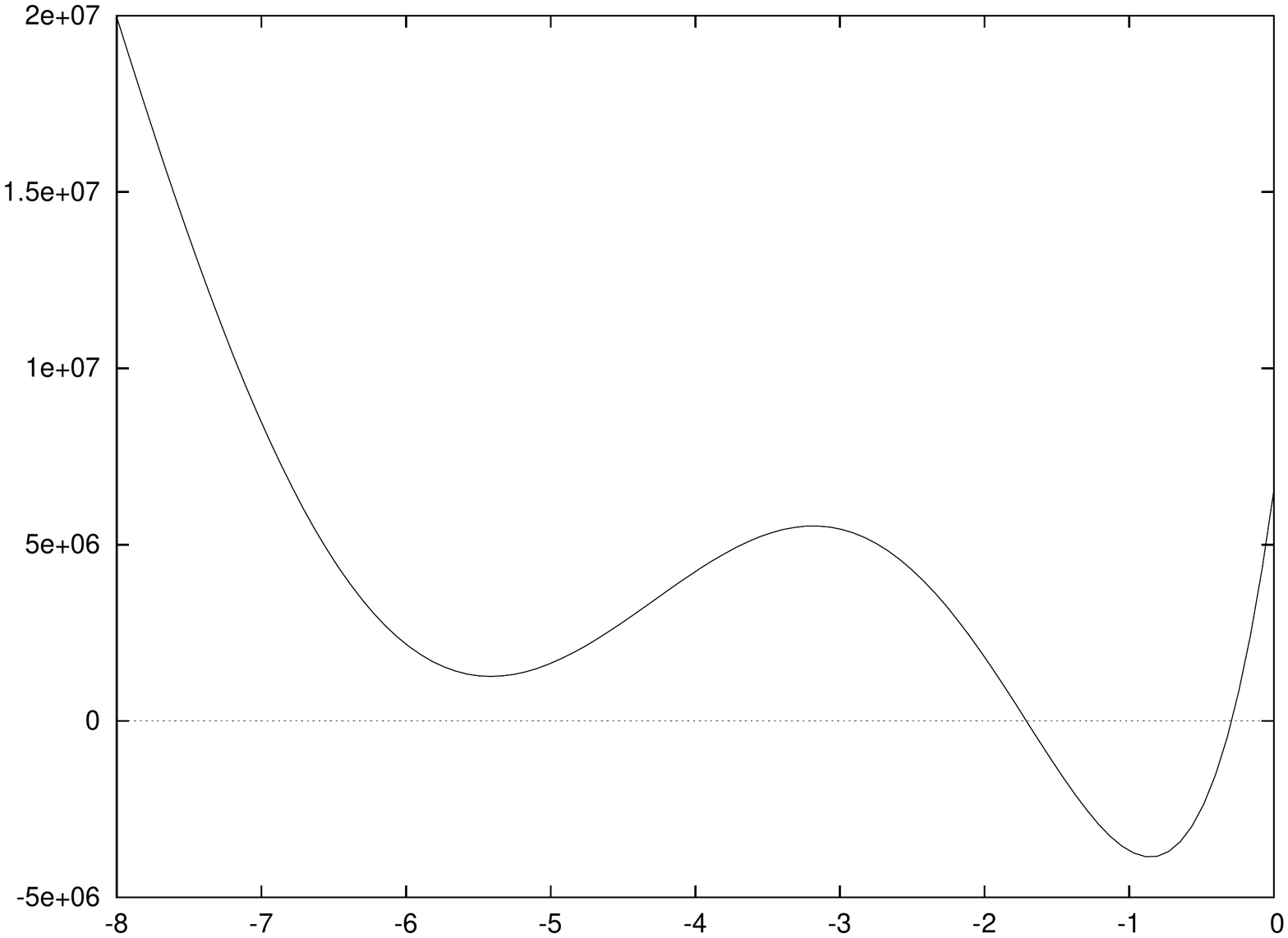}
\caption{Graph of $R\left( x\right) $
}
\end{figure}

Similarly it is possible to show for
\[
P_{11}\left( x\right) =\frac{x}{11!}\left( x+1\right) \left( x+2\right) \left( x+3\right) \left( x+8\right) \tilde{R}\left( x\right)
\]
with $\tilde{R}\left( x\right) $ an irreducible
polynomial of degree~$6$ that it has only real roots.
This can be done by checking the limits of the
intervals $\left( -67,-66\right) $, $\left( -39,-38\right) $,
$\left( -22,-21\right) $, $\left( -17,-16\right) $,
$\left( -8,-7\right) $, and $\left( -1,0\right) $.

Since $Q_n(z)= P_n(1+z)$, we have found a counterexample to Conjecture 2 (ii). This implies that also (iii) has to be revised.
Based on numerical investigations (\cite{HNR18, HNW18}) we propose the following revised version of Conjecture 2.
\newline
\newline
{\it {\bf Conjecture $2$ (New)} 
\newline
Let $n$ be a positive integer. Then the polynomial
\begin{equation}
P_n(z) = 
\sum_{ \lambda \vdash n} 
\prod_{ h \in \mathcal{H}(\lambda)} 
\left(  \frac{h^2 + z-1}{h^2} \right) \in \mathbb{C}[z]
\end{equation}
has (i) only simple roots and (ii) real part of all non-trivial roots is negative.}

\begin{remark} \ \newline
a) The first part of the conjecture has been proven for integral roots and $n$ or $n-1$ equal to a prime power \cite{HN18}.   \newline
b) The second part of the conjecture has been verified for $n\leq 700$ (\cite{HNR18}). 
Polynomials satisfying (ii) are called Hurwitz polynomials or stable polynomials. They play an important role in the theory of dynamical systems.
\end{remark}

\section{Unimodality}
Let $a_0,a_1, \ldots,a_n$ be a finite sequence of nonnegative real numbers. The sequence is denoted unimodal if
$a_0 \leq a_1 \leq \ldots \leq a_{k-1} \leq a_k \geq a_{k+1} \geq \ldots \geq  a_n$ for some $k$. It is denoted log-concave if $a_j^2 \geq a_{j-1}\, a_{j+1}$ for all $j >0$. A sequence is called ultra-log-concave if the attached sequence $a_k/ \binom{n}{k}$ is log-concave.
Due to Newton (1707), a finite sequence $a_0,a_1, \ldots,a_n$ of nonnegative real numbers with real roots is log-concave.
Actually it is already ultra-log-concave.
The $Q_n(x)$ are polynomials with nonnegative real coefficients, 
but with possible non-real roots. This makes Conjecture 3 considerably complicated.

Nevertheless, numerical calculations provide the following result.
\begin{theorem}
Let $1\leq n \leq 1000$. Then $Q_n(x)$ is ultra-log-concave. This implies Conjecture 3 (for $n\leq 1000$).
\end{theorem}
Note that for a general polynomial $P\left( x\right) $
we do \emph{not\/} have the property: $P\left( x\right) $ is unimodal if
and only if $P\left( x+1\right) $ is unimodal. For example
$x^2+2$ is not unimodal as $a_{2}=1>0=a_{1}<a_{0}=2$ but
$\left( x+1\right) ^{2}+1=x^{2}+2x+3$ is unimodal.


\begin{thebibliography}{MMM99}
\bibitem[ASW52]{ASW52}
M. Aissen, I.~J. Schoenberg, and A. Whitney: \emph{On the generating functions of totally
positive sequences.}
J. Anal. Math. \textbf{2} (1952), 93--103. 


\bibitem[Am15]{Am15} T. Amdeberhan: 
\emph{Theorems, problems and conjectures.} 
arXiv:1207.4045v6 [math.RT] 13 Jul 2015.

\bibitem[AE04]{AE04} G. E. Andrews, K. Eriksson: \emph{Integer Partitions.\/} Cambridge University Press (2004).

\bibitem[Co74]{Co74} L. Comter: \emph{Advanced combinatorics.\/} Enlarged edition,
D. Reidel Publishing Co., Dordrecht (1974).


\bibitem[DA13]{DA13} F. D'Arcais: 
\emph{D\'{e}veloppement en s\'{e}rie.} Interm\'{e}diaire Math. 
\textbf{20} (1913), 233--234.


\bibitem[Fu97]{Fu97} W. Fulton: \emph{Young Tableaux.}
Cambridge University Press 1997.



\bibitem[GKS90]{GKS90} F. Garvan, D. Kim, D. Stanton: 
\emph{Cranks and t-cores.} 
Invent. Math.
\textbf{101} (1990), 1--17.

\bibitem[Ha09]{Ha09} G.~N. Han: 
\emph{Some conjectures and open problems on partition hook lengths.} 
Exp. Math.
\textbf{18} (2009), 97--106.

\bibitem[Ha10]{Ha10} G. Han: 
\emph{The Nekrasov--Okounkov hook length formula: refinement, elementary proof and applications.} 
Ann. Inst. Fourier (Grenoble) 
\textbf{60} 
no.~1 (2010), 1--29.

\bibitem[HLN18]{HLN18} B. Heim, F. Luca, M. Neuhauser: 
\emph{On cyclotomic factors of polynomials related to modular forms.} The Ramanujan Journal (2018) DOI:
10.1007/s11139-017-9980-8.


\bibitem[HN18]{HN18} B. Heim, M. Neuhauser: \emph{Polynomials related to powers of the Dedekind eta function.} preprint, submitted.




\bibitem[HNR18]{HNR18} B. Heim, M. Neuhauser, F. Rupp: \emph{Imaginary powers of the Dedekind eta
function.\/} Experimental
Mathematics (2018) DOI: 10.1080/10586458.2018.1468288.

\bibitem[HNW18]{HNW18} B. Heim, M. Neuhauser, A. Weisse: 
\emph{Records on the vanishing of Fourier coefficients of powers of the Dedekind eta function.} Research in Number Theory (2018) DOI:
10.1007/s40993-018-0125-y.







\bibitem[KRY09]{KRY09} J. Kung, G. Rota, C. Yan: \emph{Combinatorics: The Rota way.}
Cambridge University Press (2009).




\bibitem[Le47]{Le47} D. Lehmer: \emph{The vanishing of Ramanujan's $\tau(n)$.} Duke Math.\ J.
 \textbf{14} (1947), 429--433.

\bibitem[Ma95]{Ma95} I.G. Macdonald: \emph{Symmetric functions and Hall polynomials.}
Second Edition, Clarendon Press, Oxford (1995).



\bibitem[NO06]{NO06} N. Nekrasov, A. Okounkov: \emph{Seiberg--Witten theory and random partitions.}
{\it The unity of mathematics.\/}  Progr. Math. {\bf{244}} Birkh\"{a}user Boston (2006), 525--596.

\bibitem[Ne55]{Ne55} M. Newman: \emph{An identity for the coefficients of certain modular forms\/}.
J. London Math.\ Soc.\ 
{\bf{30}} (1955), 488--493.



\bibitem[On03]{On03} K. Ono: \emph{ The Web of Modularity: Arithmetic of the Coefficients of Modular Forms and q-series.}
Conference Board of Mathematical Sciences 
{\bf{102}} (2003).



\bibitem[Se85]{Se85} 
J. Serre:
\emph{Sur la lacunarit\'{e} des puissances de $\eta $\/}.
Glasgow Math.\ J. 
{\bf{27}} (1985), 203--221.

\bibitem[St99]{St99} R. Stanley: \emph{ Enumerative  Combinatorics, vol.~2.}
Cambridge: Cambridge University Press, 1999.




\bibitem[We06]{We06} B. Westbury: \emph{Universal characters from the Macdonald identities\/}.
Adv. Math.\ 
{\bf{202}} 
no.\ 1 
(2006), 50--63.






 \end{thebibliography}
\end{document}